\documentstyle[12pt]{article}

\newcommand{\argmax}{\mathop{\rm argmax}\limits}

\newcommand{\const}{\mathop{\rm const}\limits}

\newcommand{\mod}{\mathop{\rm mod}\limits}

\newcommand{\Law}{\mathop{\rm Law}\limits}

\newcommand{\Variation}{\mathop{\rm Variation}\limits}

\newcommand{\Var}{\mathop{\rm Var}\limits}

\begin{document}

\begin{center}

{\bf  Rate of convergence in the   } \\

\vspace{4mm}

{\bf  maximum likelihood estimation   } \\

\vspace{4mm}

{\bf  for partial discrete parameter,  } \\

\vspace{4mm}

{\bf with applications to the cluster analysis and philology.} \\

\vspace{4mm}

 $ {\bf E.Ostrovsky^a, \ \ L.Sirota^b,  \ \ A.Zeldin^c.} $ \\

\vspace{4mm}

$ ^a $ Corresponding Author. Department of Mathematics and Computer science, \\
Bar-Ilan University, 84105, Ramat Gan, Israel.\\

\vspace{4mm}

E - mail: \ galo@list.ru \  eugostrovsky@list.ru\\

\vspace{4mm}

$ ^b $  Department of Mathematics and computer science. Bar-Ilan University, 84105,\\ Ramat Gan, Israel.\\

\vspace{4mm}

E - mail: \ sirota3@bezeqint.net.il\\

\vspace{4mm}

$ ^c $  Research and consulting officer, the Ministry of Immigrant Absorption, Israel. \\

\vspace{4mm}

E - mail: \ anatolyz@moia.gov.il\\

\vspace{5mm}
                    {\bf Abstract.}\\

 \end{center}

 \vspace{4mm}

  The problem of estimation of the distribution parameters on the sample when the part of these parameters are discrete
  (e.g. integer) is considered. We prove that the rate of convergence of MLE estimates under the natural conditions on
  the distribution density is exponentially fast. \par

   We describe also the possible of the applications of the estimates offered in the cluster analysis and consequently in the
    technical diagnosis, demography and especially in philology. \par

\vspace{4mm}

{\it Key words and phrases:} Maximum Likelihood Estimation (MLE), metric entropy by Kolmogorov, relative entropy of Kullback and Leibler,
Hellinger's entropy and integral, random process (field),
stable distributions, tail of distribution, heavy tail distribution,
exponential estimation  for random fields and for sums of r.v., rate of convergence, random variables and vectors  (r.v.), large deviations,
cluster and  cluster analysis, objective function, quasi-Gaussian distribution, contrast function, nuisance parameter, action function,
mixture, density of distribution, Cartesian and polar coordinates.\\

\vspace{4mm}


\vspace{4mm}

\section{Introduction. Statement of the problem. Notations. Definitions. }

\vspace{3mm}

   Let $ (\Omega, \cal{B}, {\bf P}) $ be a probability space with the expectation $ {\bf E,}$
 and $ (X, \cal{A}, \mu) $  be a measurable space with sigma-finite non-trivial measure $ \mu, $ and
$  \Theta $ be arbitrary  locally  compact topological space equipped by the ordinary Borelian
sigma - field $  G, $  whereas  $\ F = \{f \}, \ f = f(x, \theta), \ \theta \in \Theta, \ x \in X $
be a family of a {\it strictly positive}
 $ (\mod \mu) $ probabilistic densities:

$$
\int_X f(x,\theta) \ \mu(dx) = 1, \ \theta \in \Theta
$$
to be assumed continuous relative to the argument $ \theta $ for almost all values $ x; x \in X.$

 We premise also the following natural condition of the identifiability condition:

$$
\forall \theta_1, \theta_2 \in \Theta, \ \theta_1 \ne \theta_2 \Rightarrow \mu \{ x: f(x,\theta_1) \ne f(x,\theta_2)  \} > 0.
$$

 Let further $  \theta_0 $ be some fixed value of the parameter $  \theta. $ We assume that the r.v. $ \xi = \xi(\omega) $ is
a random variable (r.v) (or more generally random vector) taking the values in the space
$ X $ with the density of distribution $ f(x, \theta_0) $ relative to the measure $ \mu : $

$$
{\bf P} (\xi \in A) = \int_A f(x, \theta_0) \ \mu(dx).
$$

 Statistically this means that r.v.$ \xi $ is the (statistical) observation (or observations) with the
density $   f(x, \theta_0) $ relative to the measure $  \mu, $ where the value $ \theta_0 $ is the true, but in general case
the values of the parameter are unknown. \par

   The $  \hat{\theta} $ denotes  the Maximum Likelihood Estimation (MLE) of the parameter $ \theta $
 based on the observation $ \xi: $

 $$
 \hat{\theta} = \argmax_{\theta \in \Theta} f(\xi,\theta), \eqno(1.0)
 $$
or equally

$$
 \hat{\theta} = \argmax_{\theta \in \Theta} L(\xi,\theta), \eqno(1.1)
 $$
where the function

$$
L = L(\xi, \theta) \stackrel{def}{=} \log [f(\xi, \theta)/f(\xi, \theta_0)]
$$
is termed  the {\it contrast function,} in contradiction to the function $ \theta \to f(\xi, \theta) $  or
$ \theta \to \log f(\xi, \theta),$ which is called the ordinary Likelihood function.\par

In the  case that $  \hat{\theta} $ is not unique, we accept any arbitrary but
 measurable version of $  \hat{\theta} $   as a capacity $  \hat{\theta} $  provided it satisfies the condition (1.1.)

\vspace{3mm}

{\it  We will consider in the sequel only the case sample:   } $ \xi = \vec{\xi} = \{ \xi_1, \xi_2, \ldots, \xi_n   \}, $
{\it  where $ \xi_i, \ i = 1,2,\ldots,n $ are independent identically distributed with density    } $ f(x,\theta) $
{\it with the true value of the parameter} $ \theta = \theta_0.  $ \par

\vspace{3mm}

 It is well known that if the set $ \Theta $ is convex non-empty smooth submanifold of the whole space $  R^d, \ d = 1,2,\ldots $
and the density $  f(x,\theta) $ is in addition smooth   (of a class $  C^2 $  ) function in relation to the parameter $ \theta, $
then for all sufficiently large values $  n $ (volume of the sample) the MLE estimate $ \hat{\theta} = \hat{\theta}_n  $ based
on the whole sample does exists. This estimate  is asymptotically unbiased, asymptotically normal and is asymptotically effective with the speed
of convergence $ 1/\sqrt{n}; $ see for example \cite{Ibragimov1}. \par

 Evidently, the MLE estimation $ \hat{\theta} = \hat{\theta}_n  $ is the solution of system of equations

 $$
 \sum_{i=1}^n \frac{\partial \log \left[ f(\xi_i, \theta )/f(\xi_i,\theta_0) \right] }{ \partial \theta_k  } = 0, \ k = 1,2,\ldots,d. \eqno(1.2)
 $$

or equally
 $$
 \sum_{i=1}^n \frac{\partial \log \left[ f(\xi_i, \theta)/f(\xi_i,\theta_0) \right] }{ \partial \theta_k  } /(\theta = \hat{\theta}_n) = 0,
 \ k = 1,2,\ldots,d. \eqno(1.2a)
 $$

 The non-asymptotical estimates for the probability of $ \sqrt{n} $ deviation of a form

$$
\sup_n {\bf P} \left( \sqrt{n} ||\hat{\theta}_n - \theta_0|| > u \right) \le C_1 \exp \left(- C_2 u^{\gamma}   \right), \
u \ge 2, \ \gamma = \const > 0 \eqno(1.3)
$$
in the considered case was obtained in \cite{Ostrovsky2}. \par

\vspace{3mm}

{\bf   We consider in this article  the case when some part of the estimated parameters are discrete, (for the sake of definiteness, integer),
  and investigate the speed of convergence  MLE estimation $ \hat{\theta}_n   $  to the true value $ \theta_0. $ }\par

\vspace{3mm}

 To make the notations clearer, we accept several changes  of notations.

 $$
 \Theta = [0,1,\ldots, N] \otimes  \cal{B}, \eqno(1.4)
 $$
where $  N \in  \{1, 2,  \ldots, \infty \}, \  \cal{B} = \{  \beta \} $ is  arbitrary
separable compact topological space equipped by the ordinary Borelian sigma-field. \par

\vspace{3mm}

{\it Let the point $ \theta_0 = (0, \beta_0) $ be the true value of the parameter $ \theta, $ so that the
sample $ \xi = \vec{\xi} = \{  \xi_i  \}, \ i = 1,2, \ldots,n $ consists on the i., i.d. r.v. $ \xi_i $
 with the density } $ f = f(x,\theta_0) = f(x, 0, \beta_0). $\par

 \vspace{3mm}

  The MLE estimate  $ \hat{\theta}_n $  of the parameter $ \theta $ will  denoted by $ ( \hat{\tau}_n, \hat{\beta}_n ):  $

$$
( \hat{\tau}_n, \hat{\beta}_n ):  = \argmax_{m,\beta} \left\{ \sum_{i=1}^n \log[ f(\xi_i,m,\beta)  ] \right\}=
$$

$$
\argmax_{m,\beta} \left\{ \sum_{i=1}^n \log[ f(\xi_i,m,\beta)/f(\xi_i,0,\beta_0)   ] \right\}.\eqno(1.5)
$$

 Note that the discrete parameter estimates, without nuisance parameters, was considered in many works; see for example,
 \cite{Choirat1},  \cite{Gersanov1},  \cite{Lindsay1},  \cite{Meeden1} etc. \par

 The example of these statement of problem is described  in the articles \cite{Ostrovsky1}, \cite{Ostrovsky4}, where some problems
 of cluster analysis are considered, as well as their applications in technical diagnosis, demography and especially in philology.

 In these  cases  the number of clusters acts as a discrete parameters. \par

 Namely,  there are some grounds to accept the function $ f(x; \vec{\theta} ) $ as a density
 in technical  diagnosis, demography, philology etc.  \cite{Ostrovsky1}, \cite{Ostrovsky4}, whereas the density
 has a form of the so-called mixed quasy-Gaussian distribution (1.3)

 $$
 f(x; \vec{\theta} ) =
\sum_{k=1}^N  W_k \ G \left( x_1-a_1^{(k)}, x_2-a_2^{(k)}, \ldots, x_d-a_d^{(k)}; \vec{\alpha}^{(k)}, \ \{ \sigma_j^{(k)} \}, \ \vec{ C_1 }^{(k)} \right).
$$

 Here

 $$
\vec{\theta} = \theta \stackrel{def}{=} \{ N; \{ \vec{a}_d \}, \ \{ \vec{\sigma}_j \}, \  \{ \vec{ C_1 }^{(k)} \}   \}, \ d = \dim X, \ j,k = 0,1,2,\ldots,N,
$$

 $ W_k, \ k = 1,2, \ldots, N $ be positive numbers (weights) such that $ \sum_{k=1}^N W_k = 1. $

 \hspace{3mm}

  The quasi-Gaussian distribution was defined as follows. \par

 We denote as trivial  for any measurable set $ A, \ A \subset R $ its indicator function by $ I(A) = I_A(x):   $

$$
I_A(x) = 1, \ x \in A; \hspace{5mm} I_A(x) = 0, \ x \notin A.
$$

Let us introduce a {\it family} of  functions

$$
\omega_{\alpha}(x) = \omega_{\alpha}(x; C_1, C_2) := C_1 \ |x|^{\alpha(1)} \ I_{(-\infty,0)}(x) + C_2 \ x^{\alpha(2)} \ I_{ (0,\infty)}(x),
$$

$$
x \in R, \ C_{1,2}= \const \ge 0, \ \alpha = \vec{\alpha}  = (\alpha(1), \alpha(2)),  \ \alpha(1), \alpha(2) = \const > -1,
$$
so that $ \omega_{\alpha}(0) = 0,  $ and a family of a correspondent probability densities of a form

$$
g_{\alpha, \sigma}(x) = g_{\alpha, \sigma}(x; C_1, C_2)  \stackrel{def}{=} \omega_{\alpha}(x; C_1, C_2) \ f_{\sigma}(x),
$$

$$
f_{\sigma}(x) = (2 \pi)^{-1/2} \ \sigma^{-1} \ \exp \left(  -x^2/(2 \sigma^2)  \right).
$$

 Since

$$
I_{\alpha(k)}(\sigma) := \int_0^{\infty} x^{\alpha(k)} \exp \left( -x^2/(2 \sigma^2)  \right) \ dx = 2^{(\alpha(k) - 1)/2  } \
\sigma^{(\alpha(k) + 1)} \ \Gamma((\alpha(k) + 1)/2),
$$
where $  \Gamma(\cdot) $ is ordinary Gamma function, there is the interrelation between the constants $ C_1, C_2:  $

$$
C_1 \ I_{\alpha(1)}(\sigma) + C_2 \ I_{\alpha(2)}(\sigma) = \sigma \ (2 \pi)^{1/2},
$$
has only one degree of freedom.  In particular, the constant $  C_1 $ may be equal  to zero;
in this case the r.v. $ \xi $ possess  only non-negative values.\par

\vspace{3mm}

{\it  We will denote in the sequel  by $  C_i, K_j $ some finite non-negative constants that are not  necessary to be
the same in different places. }\\

\vspace{3mm}

 The one - dimensional distribution of a r.v. $  \xi   $ with density function of a form
$ x \to g_{\alpha, \sigma}(x - a; C_1, C_2), \ a = \const \in R $ is said  to be quasi-Gaussian or equally quasi-normal. Notation:

$$
\Law(\xi) = QN(a,\alpha,\sigma, C_1, C_2).
$$

\vspace{3mm}
Let us explain the "physical" sense of introduced parameters of these distributions.
 The value $  "a" $  may be named  {\it  quasi - center  } by analogy with normal distribution; the value $ "\alpha" $
expresses the degree of concentration of this distribution about the center  and the value of
 $ "\sigma" $ which may be  called {\it  quasi-standard } of the r.v. $  \xi  $  expressed alike in the classical Gaussian
r.v.  the degree of scattering.   \\

\vspace{2mm}

  Many properties of these distributions are previously studied in \cite{Ostrovsky1}: moments, bilateral tail behavior etc. In particular,
 it is proved that  if the r.v. $  (\xi, \eta) $ are independent and both have the quasi-Gaussian distribution
with parameters $  a = 0, \ b = 0 $ ("quasi - centered" case):

$$
\Law(\xi) = QN(0,\alpha,\sigma, C_1, C_2), \hspace{5mm} \Law(\eta) = QN(0,\beta,\sigma, C_3, C_4)
$$
may occur with different parameters $ \alpha \ne \beta, \ C_1 \ne C_3, C_2 \ne C_4  $ but with the same value of the standard
$  \sigma, \ \sigma > 0, $ then their polar coordinates  $ (\rho, \zeta)  $ are also independent. \par

  The opposite conclusion was also proved in \cite{Ostrovsky1}: the  characterization of quasi-Gaussian distribution in the demography and
 philology: if the  polar and Decart (cartesian) coordinates  are independent, then under some natural conditions
 the random variables  $ \xi, \eta $  have quasi-Gaussian distribution, and is explained why this property denotes this distribution
 of the words parameters in many languages.  \par

 \vspace{3mm}

 It is possible to  generalize our distributions  on the  multidimensional case. Actually, let us consider the  random vector
$ \xi = \vec{\xi} =  (\xi_1, \xi_2, \ldots, \xi_d) $  with the density

$$
f_{\xi}(x_1, x_2, \ldots, x_d) = G( x_1, x_2, \ldots, x_d; \vec{\alpha}, \vec{\sigma}, \vec{ C_1 },  \vec{C_2  } ) \stackrel{def}{=}
$$

$$
\prod_{j=1}^d  g_{\alpha_j, \sigma_l}(x_j; C_1^{(j)}, C_2^{(j)}),
$$
where $ \alpha_j > -1, \ \sigma_j = \const  > 0, \ C_i^{(j)}  = \const \ge 0, $

$$
C_1^{(j)} \ I_{\alpha^{(j)}(1)}(\sigma_j)  + C_2^{(j)} \ I_{\alpha^{(j)}(2)}(\sigma_j) = \sigma_j \ (2 \pi)^{1/2},
$$

 \vspace{3mm}

 An important note: during the investigation   of these discrete estimates the so-called
  Large Deviations Principle (LDP) was used \cite{Deuschel1}, \cite{Freidlin1}, \cite{Kallenberg1},  \cite{Varadhan1}; \
  \cite{Bahadur1}, \cite{Borovkov1},   \cite{Donsker1}, \cite{Dembo1},  \cite{Piterbarg1}, \cite{Puhalskii1} etc. \par

\vspace{3mm}

 We introduce some  notations. Let  $  f = f(x)  $ and  $  g = g(x) $ be two densities in relation to the measure $ \mu, $
i.e. measurable non - negative functions  such as

$$
\int_X f(x) \ d \mu = \int_X g(x) \ d \mu = 1.
$$
   A relative entropy by  Kullback and Leibler \cite{Kullback1} $  H_r(f;g) $ is defined as ordinary by the equality

 $$
 H_r(f;g) = \int_X f(x) \ \ln \left[ \frac{f(x)}{g(x)}   \right] \ \mu(dx).\eqno(1.6)
 $$
It is well known that $ H_r(f;g) \ge 0 $ and  $ H_r(f;g) = 0 $  iff $  f(x) = g(x) $  almost everywhere. \par

 Let  $  f = f(x), \  g = g(x) $  and $  h = h(x) $ be {\it three}  densities relative the measure $ \mu. $  We define a three term
relative entropy $ H_R(f; g,h) $  as follows:

  $$
 H_R(f;g,h) = \int_X f(x) \ \ln \left[ \frac{g(x)}{h(x)}   \right] \ \mu(dx). \eqno(1.7)
 $$
 Evidently, $ H_R(f;f,h) = H_R(f;h), \ H_R(f;g,h) = -H_R(f;h,g).  $\par

Entropy by Hellinger or Hellinger's integral defined for any  real number $ \lambda \in R $ and two densities $ f(x), \ g(x) $
 is by definition the following integral (if there exists)

$$
H_{ell}(\lambda; f,g) := \int_X f^{\lambda}(x) \ g^{1 - \lambda}(x) \ \mu(dx).
$$
 Hellinger   \cite{Hellinger1}  originally introduced this concept for the value $ \lambda = 1 - \lambda = 0.5. $ The general notion
 was for the first time introduced most likely in \cite{Ostrovsky2}.\par
  This notion is closely related with the so-called R'enyi and Tsallis divergences, see
  \cite{R'enyi1}, \cite{Tsallis1}. The consistent statistical estimation  of $ H_{ell}(\lambda; f,g) $ is obtained, e.g. in
  \cite{Krishnamurthy1}; see also reference therein.\par

 We offer here a slight modification of this notion, namely, a three term Hellinger's integral:

$$
H^{(3)}_{ell}(\lambda; f,g,h) := \int_X f^{\lambda}(x) \ g^{-\lambda}(x) \  h(x) \ \mu(dx), \eqno(1.8)
$$
where $  \lambda \in R, \ f,g,h $  be three densities. Of course, $  H_{ell}^{(3)}(\lambda; f,g,g) = H_{ell}(\lambda; f,g). $

 At last, let $  (T,d) $ be a metric space equipped with distance $  d = d(t,s), \ t,s \in T. $ The entropy by Kolmogorov \cite{Kolmogorov1}
 $  H(T,d, \epsilon) $ is named the natural logarithm of the minimal numbers of closed  balls $ B(y, \epsilon), \ y \in T,
 \epsilon \in (0,\infty) $ in the distance $  d  $ which cover  all the set $  T. $ \par
Obviously, $ \forall \epsilon > 0 \ \Rightarrow H(T,d, \epsilon) < \infty $  iff the set $  T $ is precompact set relative the distance $  d. $ \par

\vspace{4mm}

\section{Main result: exponential convergence for discrete parameter. }

\vspace{3mm}

 We need to  introduce some notations and conditions. $ \xi: = \xi_1,  $

$$
Q_n = Q_n(\theta) := {\bf P} (\hat{\tau}_n \ne 0) =  {\bf P} (\hat{\tau}_n \ge 1). \eqno(2.0)
$$
 This probability for confidence interval (confidence  probability) play a very important  role in our considerations.\par

  Also, let us admit

$$
a(\theta) = a(m,\beta) := {\bf E} \ln \left[ \frac{f(\xi_i,m,\beta)}{f(\xi_i,0,\beta_0)} \right] =
$$

$$
\int_X  \ln \left[ \frac{f(x,m,\beta)}{f(x,0,\beta_0)} \right] \cdot f(x,0,\beta_0) \ \mu(dx). \eqno(2.1)
$$
 The function $  -a = -a(\theta) = -a(m,\beta) $ is  relative entropy of the density  $ f(x,m,\beta) $ in relation to other density $ f(x,0,\beta_0): $

 $$
 -a(m,\beta) = H_r(f(x,0,\beta_0); f(x,m,\beta) ) =  H_r(f(x,\theta_0); f(x,\theta)).
 $$

 Therefore  $ \forall m \ge 1 \ a(m, \beta) < 0. $  Obviously, $ a(0,\beta_0) = 0.  $ \par
 We suppose in addition the function  $ \theta \to a(\theta) $ there exists and is continuous:
$  a(\cdot) \in C(\Theta_1); $  and we denote  by $  \Theta_1 $ the (closed) subspace of the space $ \Theta  $
of the form

$$
\Theta_1  = \{1,2,\ldots,N \} \otimes \cal {B}.
$$

 Moreover, we assume the following random  processes (field) (r.f.)

$$
\eta_i(m, \beta):=  \ln \left[ \frac{f(\xi_i,m,\beta)}{f(\xi_i,0,\beta_0)} \right], \
\eta^{(0)}(\theta) = \eta(\theta) - a(\theta),  \eqno(2.2)
$$

$$
 \eta(m, \beta):= \eta_1(m, \beta) = \eta(\theta), \ \theta \in \Theta_1 = \{  m, \ \beta  \}, \ m = 1,,2,\ldots,N
$$
 belong also to the space of all continuous functions $  C(\Theta_1) $ equipped  by ordinary  norm:

 $$
\forall g(\cdot) \in C(\Theta_1) \hspace{4mm} ||g|| = \max_{\theta \in \Theta_1} |g(\theta)|.
 $$

 Let $   Z = Z(A) $ be an element of conjugate space $ C^*(\Theta_1), $ i.e. countable additive signed measure defined
on the Borelian sigma-field $  G_1 $ (charge) with finite variation, which we denote by

$$
||Z|| = ||Z||C^* = ||Z||C^*(\Theta_1)  = \Variation(Z).
$$

  We postulate the finiteness of the logarithm of generating functional for the
r.f. $ \eta(m, \beta) $ for all the charges $  Z = Z(\cdot) $

$$
\Psi(Z) := \ln {\bf E} \exp \left( \int_{\Theta_1} \eta(\theta) \ Z( d \theta)  \right) < \infty. \eqno(2.3)
$$

  Moreover, we impose the classical in the theory of great deviations condition:

$$
\forall t \in (0,\infty) \Rightarrow {\bf E} \exp \left(  t ||\eta(\cdot)|| \right) < \infty. \eqno(2.4)
$$

 The following function, which  usually called action function, plays  a very important role in the theory of great deviations
is defined by the Young-Fenchel, or Legendre transform of $  \Psi(\cdot) $ in the space $ C^*(\Theta_1), $ which we will
denote also by $  \Psi^*: $

$$
I(g) := \sup_{Z \in C^*(\Theta_1)} \left\{ \int_{\Theta} g(\theta) \ Z(d \theta) - \Psi(Z)  \right\} =
$$

$$
 \sup_{Z \in C^*(\Theta_1)} \left\{ ( g(\cdot), Z(\cdot)) - \Psi(Z)  \right\}  = \Psi^*(g).  \eqno(2.5)
$$

 We denote also by $  U = C^+(\Theta_1) $ the set of all continuous functions $  \Theta \to R $  where

$$
g \in U \ (= C^+(\Theta_1)) \  \Leftrightarrow \max_{\theta \in \Theta_1} g(\theta) \ge 0.  \eqno(2.6)
$$
 Evidently, $  U  $ is closed set in the space $ C(\Theta_1),   $ and we denote  by  $  U^o $   its interior:

$$
U^o = \{g, \ \max_{\theta \in \Theta_1} g(\theta) > 0 \}.  \eqno(2.7)
$$

\vspace{4mm}

{\bf Theorem 2.1.} {\it  Let the listed above conditions are fulfilled. Suppose in addition  }

$$
\lim_{ ||Z||C^* \to \infty} \Psi(Z)/||Z||C^*  = \infty. \eqno(2.8)
$$
{\it Then}

$$
- \inf_{g \in U^o} \Psi^*(g) \le  \underline{\lim}_{n \to \infty} n^{-1} \ln Q_n \le  \overline{\lim}_{n \to \infty} n^{-1} \ln Q_n  \le
- \inf_{g \in U} \Psi^*(g). \eqno(2.9)
$$

 {\bf Proof} is  the same as in the article  of  Choirat  Ch. and  Seri R.
 \cite{Choirat1}, where   the case of the complete discrete
 parametric space $ \Theta $  is considered.  We need only to replace the finite-dimensional LDP (Large Deviation Principle) used in  \cite{Choirat1}
  by  infinite-dimensional version one,  see e.g.  \cite{Borovkov1}, \cite{Piterbarg1}, \cite{Puhalskii2}.\par

  Several details. Let us consider the partial sum denoting it as follows:

  $$
 S_n(m,\beta) = \sum_{i=1}^n \eta_i(m,\beta):
 $$

$$
Q_n = {\bf P} \left(  \sup_{m \ge 1} \sup_{\beta}  S_n(m,\beta) > 0 \right)=
 {\bf P} \left(  \sup_{m \ge 1} \sup_{\beta}  \frac{S_n(m,\beta)}{n } > 0 \right) =
$$

$$
{\bf P} \left(  S_n(\cdot, \cdot)/n   \in U^o \right). \eqno(2.10)
$$

 It is easy  to verify that all the conditions for LDP in the space $ C(\Theta_1) $ are satisfied.\\

This completes the proof of Theorem 2.1. \par

\vspace{3mm}

{\bf Corollary 2.1.}  Assume that the density  $  f(x,\theta)  $  and the space $  \Theta_1 $  are such as

$$
\inf_{g \in U^o} \Psi^*(g) =  \inf_{g \in U} \Psi^*(g). \eqno(2.11)
$$

 Then obviously

$$
  \lim_{n \to \infty} n^{-1} \ln Q_n = - \inf_{g \in U} \Psi^*(g). \eqno(2.12)
$$

\vspace{3mm}
 We will prove in the next section in particular that both the inequalities in theorem 2.1 are
in general case non-trivial.\par

\vspace{4mm}

\section{Non-asymptotical estimates. }

\vspace{3mm}

{\bf A. Lower bound.}\par

 We have:

 $$
 Q_n = {\bf P}\left(\sup_{m \ge 1} \sup_{\beta} \sum_{i=1}^n [\eta_i^{(0)}(m,\beta) +  a(m,\beta) ] > 0 \right) \ge
 $$

$$
  \sup_{m \ge 1} \sup_{\beta} {\bf P} \left( \sum_{i=1}^n [\eta_i^{(0)}(m,\beta) + a(m,\beta) ] > 0 \right) =
$$

$$
\sup_{m \ge 1} \sup_{\beta} {\bf P} \left( \sum_{i=1}^n \eta_i^{(0)}(m,\beta) > n \ |a(m,\beta)| \right)
\stackrel{def}{=} \sup_{m \ge 1} \sup_{\beta} Q_n^{(m,\beta)}, \eqno(3.1)
$$
where

$$
Q_n^{(m,\beta)}= {\bf P} \left( \sum_{i=1}^n \eta_i^{(0)}(m,\beta) > n \ |a(m,\beta)| \right).\eqno(3.2)
$$

 For the lower estimates of the variable $ Q_n^{(m,\beta)} $ we can apply the one-dimensional LDP, see for example the
 book of O.Kallenberg \cite{Kallenberg1}, p. 538 - 541.\par
 To implement this plan we will use the (generalized) Hellinger's integral (entropy) $   H_{ell}(\lambda; f,g), \ \lambda = \const. $
 We observe that the deviation function $ \Lambda(\lambda) = \Lambda(\lambda; m,\beta)  $ for the sequence $  \eta_i(m,\beta) $
is closely related to
Hellinger's integral:

$$
\exp\Lambda(\lambda; m,\beta) = {\bf E} \exp ( \lambda  \eta(m,\beta)) =
{\bf E} \exp \left( \lambda \log ( f(\xi; m,\beta)/f(\xi; 0,\beta_0  ) \right) =
$$

$$
\int_X f(x; 0,\beta_0) \cdot \exp \left( \lambda \log ( f(\xi; m,\beta)/f(\xi; 0,\beta_0) \right)  \ d \mu =
$$

$$
\int_X f^{\lambda}(x;m,\beta) \ f^{1 - \lambda}(x; 0,\beta_0) \ d \mu = H_{ell}(\lambda; f(x;0,\beta_0), f(x;m,\beta))).
$$

 Therefore as $  n \to \infty $

 $$
 \ln Q_n^{(m,\beta)} =  - n (1 + o(1))\Lambda^*(0; m,\beta)
 $$
and following

$$
\ln Q_n \ge  - n (1 + o(1)) \ \inf_{ m \ge 1, \beta  } \ \Lambda^*(0; m,\beta) \ge - 1.5 \ n \inf_{ m \ge 1, \beta  } \ \Lambda^*(0; m,\beta).
\eqno(3.3)
$$

 \vspace{3mm}

 {\bf B. Upper bound.}\par

\vspace{3mm}

 We need to introduce some new notations.

  $$
  \Delta H = H_R = H_R( m_1, \beta_1; m_2, \beta_2  ) = H_R( \theta_1;  \theta_2  )=
  $$

   $$
\Delta H(m_1,\beta_1;  m_2, \beta_2) := H_r(f(x; m_2,\beta_2); f(x; 0, \beta_0)) - H_r(f(x; m_1,\beta_1); f(x; 0, \beta_0));
 $$
 then  $ \Delta H = $

$$
\int_X f(x; 0,\beta_0) \cdot \ln \left[\frac{f(x; 0,\beta_0)}{f(x,m_1,\beta_1)} \right] \ d \mu -
\int_X f(x; 0,\beta_0) \cdot \ln \left[\frac{f(x; 0,\beta_0)}{f(x,m_2.\beta_2)} \right] \ d \mu =
$$

$$
\int_X f(x; 0,\beta_0)  \cdot \ln \left[\frac{f(x; m_2,\beta_2)}{f(x,m_1, \beta_1)} \right] \ d \mu=
H_R(f(\cdot; 0,\beta_0); f(\cdot; m_2,\beta_2), f(\cdot; m_1, \beta_1) ). \eqno(3.4)
$$

 Further, let us introduce the  following functions:

$$
\phi = \phi(\lambda; \theta_1, \theta_2), \ \overline{\phi} = \overline{\phi}(\lambda; \theta_1, \theta_2), \
$$

$$
\nu = \nu(\lambda), \ \gamma = \gamma(\lambda; \theta_1, \theta_2),  \ \overline{\gamma} = \overline{\gamma}(\lambda; \theta_1, \theta_2),
$$

$$
 \ \lambda \in R,  \ \theta_{1,2} \in \Theta_1,
$$
as follows:

$$
\phi(\lambda, \theta) := \lambda H_r + \ln \left[  H_{ell}(\lambda; f(x,m,\beta), f(x,0,\beta_0) ) \right];\eqno(3.5)
$$

$$
\overline{\phi}(\lambda, \theta) := \sup_n \left[n \phi(\lambda/\sqrt{n},\theta ) \right]; \
\nu(\lambda) := \sup_{\Theta \in \Theta_1} \overline{\phi}(\lambda, \theta); \eqno(3.6)
$$

$$
\gamma(\lambda; \theta_1, \theta_2) := \lambda H_R(\theta_1, \theta_2) + \ln H_{ell}^{(3)}(\lambda; f(x,0,\beta_0), f(x,m_1,\beta_1),
f(x,m_2,\beta_2));  \eqno(3.7)
$$

$$
\overline{\gamma}(\lambda; \theta_1, \theta_2) :=  \sup_n [n \gamma(\lambda/\sqrt{n}; \theta_1, \theta_2)].
$$

 The function $ \overline{\phi}(\lambda, \theta) $ and analogously $ \overline{\gamma}(\lambda, \theta) $ means as will be described.
If the centered random variable $ Y $  is such that

$$
{\bf E} e^{ \lambda Y } \le e^{\phi(\lambda)  },
$$
 and $ \{ Y_j \} $ are independent copies of $ Y, $ then

$$
\sup_n {\bf E} e^{\lambda \cdot n^{-1/2} \sum_{j=1}^n Y_j} \le e^{\overline{\phi} (\lambda) }.
$$

 \vspace{3mm}

 It will be assumed later that the function  $ \nu(\cdot), $
 and as well as  the function  $ \gamma(\cdot; \cdot, \cdot) $ are finite at least in some non-trivial neighborhood of origin:

$$
\exists \lambda_0 > 0 \ \Rightarrow \nu( \lambda_0) < \infty. \eqno(3.8)
$$
 We will accept

 $$
 \lambda_0 = \sup\{ z: \nu( z) < \infty \},
 $$
  as a capacity of the value $ \lambda_0  $ its maximal value; may be $ \lambda_0  = \infty. $ \par

 Let us  define a distance $ d = d(\theta_1, \theta_2) $ on the set $ \Theta_1 $ as follows:

$$
d(\theta_1, \theta_2) := \sup_{ \lambda: 0 < \lambda < \lambda_0 }
\left[  \frac{\nu^{-1}( \overline{\gamma}(\lambda; \theta_1,\theta_2 ) ) }{\lambda}  \right],\eqno(3.9)
$$
so that

$$
\overline{\gamma}(\lambda; \theta_1,  \theta_2 ) \le \nu(\lambda \cdot d(\theta_1, \theta_2)). \eqno(3.10)
$$

 Further, define

$$
G(\delta) := \sum_{m=1}^{\infty} \delta^{m-1} \ H(\Theta_1, d, \delta^m).
$$

$$
H_r :=H_r(f(x; 0, \beta_0); f(x,m,\beta));
$$

$$
\underline{H}_r := \inf_{ m \ge 1; \beta \ne \beta_0 } H_r(f(x; 0, \beta_0); f(x,m,\beta))=
 \inf_{ m \ge 1; \beta \ne \beta_0 } H_r;
$$

$$
M(u) := \inf_{\delta \in (0,1)} \left[ G(\delta)  - \gamma^*( \underline{H}_r \cdot (1 - \delta)) \right],
$$
where $ \gamma^*(\cdot) $ denotes the classical Young - Fenchel, or Legendre transform for the function $ \gamma(\cdot): $

$$
\gamma^*(u) = \sup_{ z > 0} (uz - \gamma(z)).
$$

\vspace{3mm}

{\bf  Theorem 3.1.} {\it  Assume that }  $ \underline{H}_r  > 0 $ {\it  and that } $ G(\delta) < \infty, \ \delta \in (0,1).  $
 { Then }

$$
Q_n \le \exp \{ - M( \underline{H}_r \cdot \sqrt{n} ) \}. \eqno(3.11)
$$

\vspace{3mm}

{\bf Corollary 3.1.}  Suppose in addition $  M(u) \ge K \cdot u^2 $ for all sufficiently large  values $ u; \ u \ge u_0 = \const > 0. $
 Then it follows from the assertion of theorem 3.1 the exponential non-asymptotical estimation for $  Q_n: $

$$
Q_n \le \exp \{ -  K  \underline{H}_r^2 \ n  \}, \ n \ge n_0. \eqno(3.12)
$$

\vspace{3mm}

{\bf Remark 3.1.}  The condition $ \underline{H}_r  > 0 $ is automatically satisfied if for example the set $ \Theta_1 $ is
compact set relative to the distance $ d. $ \par
 The second condition $ G(\delta) < \infty $ in turn is satisfied if the space $  \Theta_1 $ has finite dimension relative the distance $  d. $ \par

\vspace{3mm}

{\bf Proof of theorem 3.1.} Note that as before

$$
Q_n  = {\bf P} \left( \max_{\theta \in \Theta_1} \sum_{i=1}^n \eta_i(\theta) > 0  \right) =
{\bf P} \left( \max_{\theta \in \Theta_1} \left(\sum_{i=1}^n \eta^0_i(\theta) - H_r(\theta) \right) > 0  \right) \le
$$

$$
{\bf P} \left(\max_{\theta \in \Theta_1}  \sum_{i=1}^n \eta^0_i(\theta ) > n \underline{H}_r \right) =
{\bf P} \left(\max_{\theta \in \Theta_1} n^{-1/2}  \sum_{i=1}^n \eta^0_i(\theta ) > \sqrt{n} \underline{H}_r \right).
$$

 Let us introduce the centered random field (more exactly, the sequence of centered random fields)

$$
\zeta_n(\theta) = n^{-1/2}  \sum_{i=1}^n \eta^0_i(\theta ), \ \theta \in \Theta_1;
$$
 then

$$
Q_n \le {\bf  P} \left( \max_{\theta \in \Theta_1} \zeta_n(\theta) > \sqrt{n} \underline{H}_r \right).
$$

  The exact exponential bounds for tail of distribution of maximum for random fields may be found, e.g. in
\cite{Ostrovsky102}; see also \cite{Ostrovsky101}, chapter 2. We have:

$$
{\bf E} e^{ \lambda \eta^0 } = {\bf E} e^{ \lambda (\eta + H_r)} = e^{ \lambda H_r } \
{\bf E} \exp \left( \lambda \ln( f(\xi,m,\beta)/f(\xi,0,\beta_0)  )   \right) =
$$

$$
 e^{ \lambda H_r } \ \int_X f(x,0,\beta_0) \cdot f^{\lambda} (x,m.\beta) \ f^{-\lambda}(x,0, \beta_0) \ \mu(dx) =
$$

$$
 e^{ \lambda H_r } \ H_{ell}(\lambda; f(x,m,\beta), f(x,0,\beta_0) ) = e^{\phi(\lambda, \theta)}.
$$

 Therefore,

$$
{\bf E} e^{\lambda \zeta_n(\theta) } \le e^{n \phi( \lambda/\sqrt{n}, \theta) } \le
 e^{ \sup_n [n \phi( \lambda/\sqrt{n}, \theta)] } = e^{ \overline{\phi}(\lambda,\theta) } \le e^{\nu(\lambda)}. \eqno(3.13)
$$

 Let us estimate  the exponential moment for   the difference

$$
\Delta \eta = \Delta \eta (\theta_1, \theta_2) = \eta(\theta_1) - \eta( \theta_2) = H_R(f(x;\theta_0); f(x;\theta_2), f(x; \theta_1)) +
\ln \left[ \frac{f(\xi; \theta_1)}{f(\xi; \theta_2) } \right].
$$

 Thus

$$
{\bf E} e^{\lambda \Delta \eta } = e^{ \lambda H_R } \ H_{ell}^{(3)}(\lambda; \ f( x; \theta_0 ), f(x;\theta_1), f(x;\theta_2) ) =
e^{\gamma(\lambda; \theta_1, \theta_2) }.
$$

 Following,

$$
{\bf E} e^{\lambda [\zeta_n(\theta_1) - \zeta_n(\theta_2)]} \le e^{ \overline{\gamma} (\lambda, \theta_1, \theta_2)}.
$$

 It follows immediately from the direct definition  of the $ d \ - $ distance that

$$
{\bf E} e^{\lambda [\zeta_n(\theta_1) - \zeta_n(\theta_2)]} \le e^{ \nu (\lambda \cdot d( \theta_1, \theta_2))}. \eqno(3.14)
$$

 The inequalities (3.13) and (3.14) may be rewritten on the language $ B(\phi) $ spaces, see
\cite{Kozatchenko1}, \cite{Ostrovsky101}, chapter 1, as follows

$$
\sup_{\theta \in \Theta_1} \sup_n ||\zeta_n(\theta)||B(\nu) \le 1, \eqno(3.15a)
$$

$$
\sup_n ||\zeta_n(\theta_1) - \zeta_n(\theta_2)||B(\nu) \le d(\theta_1, \theta_2). \eqno(3.15b)
$$

   It remains to apply the main result of \cite{Ostrovsky101}, chapter 3,  section 3.4. \par

\vspace{4mm}

\section{Examples.}

\vspace{3mm}

 {\bf A.  Regular case.} If $ \xi_i $ are Gaussian distributed with parameters $ m = {\bf E} \xi_i = 0,1,\ldots $
 and $ \beta = \Var \xi_i \ge 1 $ we conclude:

 $$
  Q_n  = 1 - \Phi( \sqrt{n}), \ \Phi(z) = (2 \pi)^{-1/2} \int_{-\infty}^z \exp (-y^2/2) \ dy.
 $$

As a consequence: there holds for suitable greatest values $  n,  $ for instance, $ n \ge 4  $

$$
(2 \pi)^{-1/2} \ n^{-1/2}  \ e^{-n/2} \ \left( 1 - \frac{c_-}{n} \right) \le Q_n  \le
(2 \pi)^{-1/2} \ n^{-1/2}  \ e^{-n/2} \ \left( 1 + \frac{c_+}{n} \right).
$$
holds true  for suitable greatest values $  n,  $ for instance, $ n \ge 4  $

\vspace{4mm}

 We will prove further that if the conditions of theorem 2.1 are not satisfied, the speed of convergence $  Q_n \to 0 $  may
differ from the exponential. \par
 Obviously, if for all the values $  m = 1,2,\ldots N \ $
 $  f_0(x) \ne f_m(x) $ on the set of positive measure, then $  \lim_{n \to \infty} Q_n = 0. $\par

\vspace{3mm}

{\bf B. Stretched exponential random variables. }

\vspace{3mm}

 The  distribution of a r.v. $  \xi $ for which

$$
c_1(t) \exp \left(  - b(t) t^r  \right) \le  {\bf P}( \xi > t  ) \le c_2(t) \exp \left(  - b(t) t^r  \right),
$$

$$
c_3(t) \exp \left(  - b(t) t^r  \right) \le  {\bf P}( \xi < - t  ) \le c_4(t) \exp \left(  - b(t) t^r  \right), \eqno(4.1)
$$

$$
 \ t > t_0 = \const > 0,  \ r = \const \in (0,1),
$$
where $  c_k(t), \ b(t) $ are  positive continuous slowly varying functions, is named in the article \cite{Gantert1}
{\it  stretched exponential  distribution. }\par

 Let for definiteness $  c_k(t) = \const > 0,  b(t) =1. $  Denote $  m = {\bf E} \xi, \ \xi_i $ be independent copies $ \xi. $ It is
proved in particular in  \cite{Gantert1} that

$$
\lim_{n \to \infty} \frac{1}{n^r} \ln {\bf P}  \left( \sum_{i=1}^n \xi_i > x  \right) = - (x-m)^r, \ x > m. \eqno(4.2)
$$
See also earlier publication of Nagaev S.V. \cite{Nagaev1}. \par

 Assume in addition that the r.v. $ \xi $ has a positive even density $  f_0 = f_0(x) = f_0(|x|), $  so that $ \xi_i $ are symmetrically
distributed and hence $  m = 0. $ \par
 Introduce a second density  $  f_1 = f_1(x) = f_1(|x|) $ as follows:

$$
f_1(x) = C \ e^{-|x|} \ f_0(x), \hspace{5mm} \int_R f_1(x) \ dx = 1,
$$
and consider the following estimation problem $ \Theta = \{0, 1 \}, \ \Theta_1 = \{ 1 \}, $ in other words, testing of statistical
hypotheses.  It follows from the cited main result  of \cite{Gantert1}  that

$$
\ln Q_n  = - C_0 n^{r} (1 - o(1)), \ n \to \infty, \ C_0 = \const \in (0,\infty). \eqno(4.3)
$$

 Note that in the article  \cite{Gantert1} is considered the case of {\it  weighted sums } of independent random variables; see also
\cite{Kiesel1}.\par

\vspace{3mm}

{\bf C. Random variables  with heavy/power tails. }

\vspace{3mm}

 We consider in this  subsection the case when $ \xi, \ \{\xi_k \} $  have symmetrical (even) density and are i.i.d.
 so as for some  $ p = \const > 2 \ \Rightarrow  $

 $$
 |\xi|_p^p \stackrel{def}{=} {\bf E} |\xi|^p < \infty.
 $$

 For instance,

 $$
 f_0(x) = \frac{C_0(p)}{(1 + |x|^{p+1} ) \ (\ln (e + |x| )^2}, \hspace{5mm} \int_R f_0(x) \ dx = 1.
 $$

  We denote also $  a = {\bf E} |\xi_k| =\const < \infty. $ \par
 Let, as before, in any case the alternative density $  f_1 = f_1(x) = f_1(|x|) $ looks like

$$
f_1(x) = C \ e^{-|x|} \ f_0(x), \hspace{5mm} \int_R f_1(x) \ dx = 1.
$$

 Then

 $$
 Q_n = {\bf P} \left( \sum_{k=1}^n \ln \frac{f_1(\xi_k)}{f_0(\xi_k)} > 0   \right) = {\bf P} \left( \sum_{k=1}^n (b - |\xi_k|) > 0  \right) =
 $$

 $$
 {\bf P} \left( \sum_{k=1}^n ( |\xi_k|  - a) < - nd   \right), \ d = \const > 0.
 $$
  We deduce the following inequality using the Rosenthal's and Tchebychev's inequalities

 $$
 Q_n \le C_R^p \frac{|\xi|_p^p \ p^p}{n^{p/2} \ d^p \ \ln^p p}, \eqno(4.4)
 $$
 where $ C_R \approx 1.773682 $ is the absolute known constant, see \cite{Ostrovsky6}.\par

  The more exact estimate for  $  Q_n $ in the case when $  p > 4  $ may be obtained from the famous theorem of Baum and Katz \cite{Baum1}
 under at the same condition $ |\xi|_p < \infty:  $

$$
Q_n \le \frac{C_3(p, |\xi|_p)}{n^{p - 2}},
$$
 but with an unknown constant $ C_3 = C_3(p, |\xi|_p).   $ \par

\vspace{3mm}

{\bf D. Random variables from Grand Lebesgue spaces. }

\vspace{3mm}

 Let $  \psi = \psi(p), \ p \in (a,b),  \ a = \const \ge 2, \ b = \const \in (a, \infty) $ be  continuous positive function
so that the function $  p \to p \cdot \log \psi(p) $  would be convex.\par

 In the case when $ b = \infty $ we impose on the function $  \psi(\cdot) $ in addition the restriction

$$
\lim_{p \to \infty} \psi(p) = \infty.
$$
 For instance, $  \psi(p) = p^l, \ l = \const > 0 $ or $ \psi(p) = p^{l_1} \ \log^{l_2} (p), l_j = \const, \ l_1 > 0. $\par

 \vspace{3mm}

 The Banach space $  G \psi $ consists by definition on all the r.v. $  \{\eta \} $ defined on the fixed probability space
with finite norm

$$
|| \eta ||G \psi := \sup_{ p \in (a,b)} \left[ \frac{|\eta|_p}{\psi(p)}  \right]. \eqno(4.5)
$$

 This spaces were introduced in  \cite{Kozatchenko1}, more detailed  investigation of these spaces may be found in the
 monograph  \cite{Ostrovsky101}, chapters 1,2.\par

  We define for all such a function $ \psi(\cdot) $ a new functions

$$
\psi_1(p) = C_R \ d^{-1} \psi(p) \ p /\ln p, \ p \in (a,b);
$$

$$
\psi_2(p) = p \ \ln \psi_1(p); \hspace{5mm} \psi_3(p) = \psi_2^*(z) = \sup_{p \in (a,b)}(p z - \psi_2(p)). \eqno(4.6)
$$
 The operator $  \psi \to \psi^* $ is called the  Young-Fenchel, (it is finite, in general case)
  or Legendre transform.\par

 We suppose as a continuation of the subsection {\bf C } that the random variables $  \xi, \ \xi_k  $ belong to some
$  G \psi  $  space. We can, for instance, choose this function $ \psi(\cdot) $  by a so-called natural way:

$$
\psi(p) := |\xi|_p,
$$
if, of course, $  |\xi|_p < \infty $ for any value $ p, \ p > 2. $ As a capacity of the boundaries $ a,b $ we put

$$
a := 2, \hspace{6mm}  b := \sup \{ p: \ \psi_p < \infty \},
$$
may be $ b = \infty. $ Naturally, in this case $ ||\xi||G\psi = 1.  $\par

 We assert if $ \xi \in G(\psi): $

$$
Q_n \le \exp \left\{ - \psi_3[ \ln (  \sqrt{n} /||\xi||G \psi ) ] \right\},  \eqno(4.7)
$$
the so-called {\it subexponential } estimate. \par

 Indeed, it follows from the direct definition of the $ G\psi $ norm

$$
|\xi|_p \le \psi(p) \ ||\xi||G\psi, \ p \in (a,b).
$$
 We deduce after substituting into (4.4)

$$
Q_n \le \frac{\psi_2^p(p) \ ||\xi||^p G\psi}{n^{p/2} } = \exp \{- [p \ln ( \sqrt{n}/||\xi||G\psi  ) - p \ln \psi_2(p)]  \}.
$$

 It remains to take the minimum over $ p. $ \par

\vspace{3mm}

{\bf E. Stable distributed   variables. }

\vspace{3mm}

 Let now $  \xi, \ \xi_k, \ k = 1,2,\ldots,n  $ be i.i.d. random variables with the density $  f_0 = f_0(x) = f_0(|x|) $
 with symmetric {\it  stable} distribution:

$$
{\bf E} e^{i t \xi_k} = \int_R e^{i t x} \ f_0(x) \ dx = e^{ -|t|^{\alpha} }, \ t \in R, \ \alpha = \const \in (1,2).
$$

The condition $  \alpha > 1 $ guarantees  the finiteness of the first moment:   $  {\bf E} |\xi| < \infty. $\par

Let us introduce as before a second density  $  f_1 = f_1(x) = f_1(|x|) $ as follows:

$$
f_1(x) = C \ e^{-|x|} \ f_0(x), \ \int_R f_1(x) \ dx = 1,
$$
and consider again the following estimation problem $ \Theta = \{0, 1 \}, \ \Theta_1 = \{ 1 \}, $ to put it differently, testing of statistical
hypotheses.  It follows from the  main result  of articles \cite{Amosova1}, \cite{Heyde1}, \cite{NagaevA1}, \cite{Nagaev1},  \cite{Pinelis1}  that

$$
 Q_n  =  \frac{C_{\alpha}}{ n^{\alpha - 1}} (1 + o(1)), \ n \to \infty, \ C_{\alpha} = \const \in (0,\infty). \eqno(4.8)
$$

 The  case of non-symmetrical distribution  that will be treated later was  investigated in \cite{Kim1}, \cite{NagaevA1}. \par

 It is interesting to note by our opinion  that when the true distribution of the sample has a density $  f_1(x), $ then
the error probability $  \tilde{Q}_n $  has as ordinary an exponential form:

$$
 |\ln \tilde{Q}_n | \asymp C_1(\alpha) \ n, \ n \to \infty.
$$

  Indeed, we have

$$
\sum_{k=1}^n \ln \left[ \frac{f_1(\xi_k)}{f_0(\xi_k)}  \right] = n C_2 - \sum_{k=1}^n |\xi_k|.
$$

  But it is known from the cited articles that for the random sequences from the domain of stable attraction
under considered condition    $  \alpha > 1 $

 $$
 {\bf P}(S_n - {\bf E} S_n  > x) \sim n {\bf P}(\xi - {\bf E }\xi > x), \ x \asymp n.
 $$

\vspace{3mm}

 In the more general case when a symmetric distributed i., i.d.  random variables $  \{ \xi_k \} $ having a heavy regular varying tail of
distribution:

$$
{\bf P} (\xi_k > x) \sim x^{-\alpha} \ L(x), \ x \to \infty,  \ \alpha = \const \in (1,2),
$$
where $  L = L(x) $ is positive continuous slowly  varying function as $  x \to \infty, $ we can conclude

$$
\lim_{n \to \infty} \frac{ {\bf P} (S_n > x + a )}{n^{1-\alpha} L(n)} = x^{-\alpha}, \ a = {\bf E} \xi, \ x \asymp n.
$$

 Therefore we deduce in the considered case

$$
Q_n \sim C_3(\alpha) \ n^{1 - \alpha} \ L(n). \eqno(4.9)
$$

\vspace{3mm}

Let us consider now  the case $ \alpha \in (0,1). $ More exactly, let the r.v. $ \xi $  obeys a standard symmetric stable distribution
with such a value of the parameter $  \alpha. $  We derive analogously using the particular case of the results of
 Amosova  \cite{Amosova1}

$$
Q_n \sim \frac{C(\alpha)}{n^{ 1/\alpha -1 }}, \ n \to \infty. \eqno(4.10)
$$

\vspace{3mm}

 At last, in symmetrical case with $ \alpha = 1, $ i.e. when $ \xi $ has a classical Cauchy distribution

$$
f_0(x) = \frac{\pi^{-1}}{1+x^2}, \ x \in (-\infty, + \infty)
$$
 and as before $ f_1(x) = C e^{-|x|} f_0(x), $ then

$$
Q_n \sim \frac{C}{\ln n}, \ n \to \infty, \ n \ge 3. \eqno(4.11)
$$

\vspace{3mm}

{\bf F. Martingale generalization. }

\vspace{3mm}

 Let us assume again $  f_1(x) = C e^{-|x|} f_0(x).  $ We continue to accept $ \theta = 0 $ if

 $$
 \sum_{i=1}^n \ln \frac{f_0(\xi_i)}{f_1(\xi_i)} > 0
 $$
and $ \theta = 1  $ otherwise. \par

  But we suppose in this subsection that the sequence  of the random variables \\
  $ \{ \eta_i \} := \{ |\xi_i| - {\bf E} |\xi_i| \}, \ i = 1,2,\ldots $ forms
 the sequence of centered martingale - differences under certain filtration $  \{  F_i \}. $ \par

 As before

$$
Q_n = {\bf P} ( \sum_{i=1}^n \eta_i > n d ), \ d = \const > 0.
$$

 The exact non-asymptotic estimations for these probabilities for martingales can be found in
 \cite{Ostrovsky5}; see also  \cite{Lesign1},  \cite{Li1}. \par

  For example, if for some $   p \ge 2 \ \forall i \ \Rightarrow |\eta_i|_p < \infty, $ then

$$
Q_n \le d^{-p} \ (p-1)^p \ n^{-p/2} \ \left\{ n^{-1} \sum_{i=1}^n |\eta_i|_p^2 \right\}^{p/2}. \eqno(4.12)
$$

 Another example. Introduce the tail function $  T = T(x), \ x> 0 $  for the sequence  $ \{  \eta_i \} $  as follows:

$$
T(x) := \sup_i \max( {\bf P} (\eta_i > x), \ {\bf P} (\eta_i < -x)  )
$$
and define

$$
W[T](x) = \min \left(1, \inf_{v > 0} \left[ e^{ -x^2/(8 v^2) } - \int_v x^2 d T(x)  \right]  \right),
$$
if of course $ \int_0^{\infty} x^2 | dT(x)| < \infty. $ \par
 Proposition:

 $$
 Q_n \le W[T](d \sqrt{n}). \eqno(4.13)
 $$

 A particular case for some $ q = \const > 0, \ K = \const > 0 $

$$
T(x) \le \exp \left( - (x/K)^q  \right), \ x > 0,
$$
then

$$
Q_n \le \exp  \left(  - C(d) n^{q/(q+2)} \ K^{-2q/(q+2)}  \right), \eqno(4.14)
$$
and the last estimate is unimprovable.\par

\vspace{3mm}

\section{ Concluding remarks.}

\vspace{3mm}

{\bf Confidence region  for "continuous" parameters.}\par

 Suppose the set $ B = \{\beta \}  $  is compact  smooth  $ (C^3)  $ subset  of the Euclidean space $  R^d $ equipped with
ordinary norm $ |\beta|. $ We consider the confidence probability

$$
W_n = {\bf P} (\sqrt{n} | \hat{\beta}_n - \beta_0| > u  ), \ u = \const > 0, \eqno(5.1)
$$
where as before $ ( \hat{\tau}_n, \hat{\beta}_n) $ is MLE estimation for $ (m, \beta). $\par

 We get:

$$
W_n = {\bf P} \left(  \hat{\tau}_n = 0, \  \sqrt{n} | \hat{\beta}_n - \beta_0| > u  \right)  +
$$

$$
{\bf P} \left(  \hat{\tau}_n \ge 1, \  \sqrt{n} | \hat{\beta}_n - \beta_0| > u  \right) \stackrel{def}{=} V_1(n) + V_2(n,u).
$$

 If all the conditions of theorem 2.1 are satisfied, then

 $$
 V_2(n,u) \le {\bf P} \left(  \hat{\tau}_n \ge 1 \right) \le e^{ - C \ n }.\eqno(5.2)
 $$

As for the probability $ V_2(n,u),  $ that

 $$
 V_2(n,u) \le
 {\bf P} \left( \max_{ |\beta - \beta_0|> u/\sqrt{n} } \sum_{i=1}^n \left[ \ln \frac{f(\xi_i,0,\beta)}{f(\xi_i,0,\beta_0)} \right] > 0  \right) \le
 $$

 $$
 \le e^{ - C_2 \ u^{\kappa} }, \ \kappa = \const \in (0, 2], \ C_2 = \const > 0, \eqno(5.3)
 $$
see \cite{Ostrovsky2}.\par

Eventually,

$$
W_n \le e^{ - C \ n } + e^{ - C_2 \ u^{\kappa} }, \ \kappa = \const \in (0, 2], \ C, C_2 = \const > 0. \eqno(5.4)
$$


\vspace{4mm}

\begin{thebibliography}{99}

\vspace{3mm}


\bibitem{Ostrovsky1}
{\sc Ostrovsky E., Sirota L., and  Zeldin A.} {\it Characterization  of quasi-Gaussian distributions.}
arXiv:1311.2341v1 [math.ST] 11 Nov 2013

\bibitem{Ostrovsky2}
{\sc Ostrovsky E., Rogover E.} {\it Non - asymptotic exponential bounds for MLE deviation under minimal conditions
via classical and generic chaining methods.}
arXiv:0903.4062v1 [math.PR] 24 Mar 2009

\bibitem{Ostrovsky4}
{\sc Ostrovsky E., Sirota L., and Zeldin A.} {\it Parametric density-based optimization of partition in cluster analysis,
with applications.}
arXiv:1312.3038v1 [math.ST] 11 Dec 2013


\bibitem{Ostrovsky5}
{\sc Ostrovsky E., Sirota L.} {\it  Tail estimates for martingale under "LLN" norming.}
arXiv:1207.1908v1 [math.PR] 8 Jul 2012


\bibitem{Ostrovsky6}
{\sc Ostrovsky E., Sirota L.} {\it Schl\"omilch and Bell series for Bessel's functions, with probabilistic applications.   }
arXiv:0804.0089v1 [math.CV] 1 Apr 2008

\vspace{10mm}

\bibitem{Anderberg1}
{\sc Anderberg, M.R. } {\it Cluster Analysis for Applications. } Academic Press, New York, 1973.

\bibitem{Hellinger1}
{\sc Hellinger, E.} (1909),
{\it Neue Begr\"undung der Theorie quadratischer Formen von unendlichvielen Ver\"anderlichen.}
Journal für die reine und angewandte Mathematik (in German), {\bf 136,}  210 \ - \ 271.


\bibitem{Kalling1}
{\sc Kailing K.,  Kriegel H.-P., and Kr\"oger P.}{\it Density-Connected Subspace Clustering for High-Dimensional Data.}
 In: Proc. SIAM Int. Conf. on Data Mining (SDM'04), pp. 246-257, 2004.

\bibitem{Kriegel1}
{\sc Kriegel H.-P., Kr\"ogel P., Sander L., and Zimek A. }
 {\it Density-based Clustering.} WIREs Data Mining and Knowledge Discovery,  (2011), {\bf 1 (3)}, 231–240.

\bibitem{Pujol1}
{\sc  Pujol J.M.,  Javier Bejar J., and  Delgado J.}
{\it Clustering algorithm for determining community structure in large networks.}
Physical Review,  E 74, (2012), 016107, (2006), 47 - 54.

\vspace{10mm}

\bibitem{Amosova1}
{\sc Amosova N.N.} {\it  Probabilities of large deviations in the case of stable  limit distribution . }
Mat. Zametki, {\bf 35,} (1984), 125 \ - \ 131.


\bibitem{Bahadur1}
{\sc Bahadur, R., Zabell, S., and Gupta, J.} (1980). {\it Large deviations, tests, and estimates.}
 In I.M. Chaterabarli (ed.), Asymptotic Theory of Statistical Tests and Estimation, pp. 33 \ - \ 64. New York: Academic Press.
Mathematical Reviews (MathSciNet): MR571334


\bibitem{Baram1}
{\sc Baram Y. and Sandell, N. R. Jr.} (1978). {\it Consistent
estimation on finite parameter sets with application to linear systems identification.} IEEE Trans.
Automat. Control, {\bf 23},  451 \ – \ 454. MR0496912

\bibitem{Baum1}
{\sc Baum L.E. and Katz M.} {\it  Convergence rates  in the law of large numbers. } Trans. of AMS, {\bf 120},
(1965), 108 \ - \ 123.


\bibitem{Borovkov1}
{\sc Borovkov A.A. and Mogul'skii A.A.} {\it  Probabilities of large deviations in topological spaces, I, II.  }
Sibirsk. Math. Zh., {\bf 19,} (1978), 988 \ - \ 1004; {\bf 21 (5),} (1980), 12 \ - \ 26;  English transl. in Siberian  Math. J.
{\bf 19,} (1978), {\bf 21,} (1980).


\bibitem{Choirat1}
{\sc Choirat  Ch. and  Seri R.}
{\it Estimation in Discrete Parameter Models.}
Statistical Science, (2012,) Vol. 27, No. 2, $ 278 \ – \ 293. $
DOI: 10.1214/11-STS371.\\
 Also: arXiv:1207.5653v1 [stat.ME] 24 Jul 2012.

\bibitem{Deuschel1}
{\sc Deuschel, J. and Stroock, D.} (1989). {\it Large Deviations.}  Boston: Academic Press.
Mathematical Reviews (MathSciNet): MR997938


\bibitem{Donsker1}
{\sc Donsker M.D.and  Varadhan S.R.S.}  {\it  Asymptotic evaluation  of certain  Markov processes  expectations for large time, I,II,III. }
Comm. Pure Appl. Math. {\bf 28,} (1975),  1 \ - \  47; {\bf 28,} (1975),  279 \ - \ 301; {\bf29,} (1976),  389 \ - \ 461.

\bibitem{Dembo1}
{\sc Dembo A. and Zeitouni O.} {\it  Large deviations techniques and applications.} Jones and Bartlett, Boston, (1993).

\bibitem{Freidlin1}
{\sc Freidlin, M.I. and Wentzell, A. D.}
(1979) {\it Random Perturbations of Dynamical Systems.} Moscow: Nauka [in Russian]. English translation: Springer (1984).
Mathematical Reviews (MathSciNet): MR722136


\bibitem{Gantert1}
{\sc  Gantert N., Ramanan K. and Rembart F.}  {\it Large deviations for weighted sums of stretched  exponential random variables.}
arXiv:1401.4577v1 [math.PR] 18 Jan 2014

\bibitem{Gersanov1}
{\sc Gersanov A. M.} (1979). {\it Optimal estimation of a discrete parameter.} Teor. Veroyatnost. i Primenen. {\bf 24},
220 \ – \ 224. MR0522259

\bibitem{Heyde1}
{\sc Heyde C.C.} {\it On large deviation problem for sums of random variables which are not attracted to the normal law. }
 Ann. Math. Stat., {\bf 38,}  (1967), 1575 \ - \ 1578.

\bibitem{Ibragimov1}
{\sc Ibragimov I.A.  and Khasminskii R.Z.} (1981). {\it Statistical estimation: Asymptotic Theory.} Springer Verlag
(Russian ed. 1979).

\bibitem{Kallenberg1}
{\sc  Kallenberg Olav.} {\it  Foundation of Probability.} (2001), Second edition,  Springer Verlag,  New York - ... - Tokyo.

\bibitem{Kiesel1}
{\sc Kiesel, R. and Stadtm\"uller, U.} (2000). {\it A large deviation principle for weighted sums of independent and
identically distributed random variables. } Journal of Mathematical Analysis, 251:929–939.

\bibitem{Kim1}
{\sc  Kim L.V. and  Nagaev A.V.} {\it The nonsymmetric problem of large deviations} (in Russian),
Teor. Veroyatnost. i Primenen. 20 (1) (1975), pp. 58 \ - \ 68.

\bibitem{Kluppelberg1}
{\sc Kl\"uppelberg  C., Mikosch T.} {\it Large Deviations on Heavy - Tailed Random Sums  with applications in Insurance and Finance.}
 J. Appl. Probab., {\bf 34,} (1997),  293 \ - \ 308.

\bibitem{Kolmogorov1}
{\sc Kolmogorov, A. N. and Tikhomirov, V. M.}  (1959), {\it $ \epsilon $ - entropy and $ \epsilon $ - capacity
of sets in a functional space. } Uspekhi Mat. Nauk, 14, 3; {\bf 86}.


\bibitem{Kozatchenko1}
{\sc  Kozatchenko Yu. V. and  Ostrovsky E.I. } {\it Banach spaces of random
variables of subgaussian type.} Theory Probab. Math. Stat., Kiev,
1985,  42--56 (Russian).

\bibitem{Krishnamurthy1}
{\sc Akshay Krishnamurthy, Kirthevasan Kandasamy, Barnab'as P'oczos  and  Larry Wasserman.  }
{\it Nonparametric Estimation of R'enyi Divergence and Friends.}
arXiv:1402.2966v1 [stat.ML] 12 Feb 2014



\bibitem{Kullback1}
{\sc Kullback, S.; Leibler, R.A. } (1951). {\it On Information and Sufficiency.} Annals of Mathematical Statistics, {\bf 22,} (1); 79 \ - \ 86.

\bibitem{Lesign1}
{\sc Lesign E., Volny D.} {\it  Large deviations for martingales. }  Stochastic Processes and
their Applications, 96, 143 - 159 (2001).

\bibitem{Li1}
{\sc Li Y. } (2003). {\it A martingale inequality and large deviations.} Statist. Probab. Lett.
{\bf 62,} 317 \ - \ 321.

\bibitem{Lindsay1}
{\sc Lindsay, B. G. and Roeder, K.} (1987). {\it A unified treatment of integer parameter models.} J. Amer. Statist.
Assoc., {\bf 82}, 758 \ – \ 764. MR0909980

\bibitem{Meeden1}
{\sc Meeden, G. and Ghosh, M.} (1981).  {\it Admissibility in finite problems. } Ann. Statist., {\bf 9,} 846 \ – \ 852. MR0619287

\bibitem{NagaevA1}
{\sc Nagaev A.V.} {\it On the asymmetric problem of large deviations when the limit law is  stable. } Theor. Probab. Appl., {\bf 28,}
(1983), 670 \ - \ 680.

\bibitem{Nagaev1}
{\sc Nagaev, S. V.} (1969). {\it Integral limit theorems for large deviations when Cram'er’s condition is not fulfilled.}
Theory of Probability and its Applications, {\bf 14}, (1), 51 \ – \ 64.

\bibitem{Nagaev2}
{\sc Nagaev, S. V.}  (1979). {\it Large deviations for sums of independent random variables.}  Annals of Probability,
{\bf 7,} 745 \ – \ 789.

\bibitem{Ostrovsky101}
{\sc  Ostrovsky E.I.} {\it Exponential Estimations for Random Fields.}
Moscow - Obninsk, OINPE, 1999 (in Russian).

\bibitem{Ostrovsky102}
{\sc Ostrovsky E.I.} (2002). {\it Exact exponential estimations for random field maximum
distribution.} Theory Probab. Appl. 45 v.3, 281 - 286.

\bibitem{Pinelis1}
{\sc Pinelis I.} (1985) {\it On the asymptotic equivalence of probabilities of large deviations for
sums and maxima of independent random variables,} (in Russian). In: Limit Theorems in Probability Theory ,
Trudy Inst. Math., {\bf 5,}  Nauka, Novosibirsk.

\bibitem{Piterbarg1}
{\sc Piterbarg V.I.  and Fatalov V.R.} {\it  The Laplace method for probability measures in Banach spaces.}
Russian Math. Surveys, 1995, {\bf  50,} 1151 \ - \ 1239.

\bibitem{Puhalskii1}
{\sc Puhalskii A. and  Spokoiny V.} {\it On large-deviation efficiency in statistical inference.}
 Bernoulli, Volume 4, Number 2 (1998), 203-272.

\bibitem{Puhalskii2}
{\sc Puhalskii, A. } (1991) {\it On functional principle of large deviations.} In V. Sazonov and T. Shervashidze (eds),
New Trends in Probability and Statistics, Vol. 1, pp. 198-218. Utrecht: VSP/Moks'las.
Mathematical Reviews (MathSciNet): MR1200917

\bibitem{Puhalskii3}
{\sc Puhalskii, A.} (1993)  {\it On the theory of large deviations.} Theory Probab. Appl., 38(3), 490-497.
Mathematical Reviews (MathSciNet): MR1404664

\bibitem{R'enyi1}
{\sc Alfr'ed R'enyi A.} {\it On measures of entropy and information.}
Proceedings of the fourth Berkeley Symposium on Mathematics, Statistics and Probability, 1960. pp. 547 \ - \ 561.

\bibitem{Rozovskii0}
{\sc Rozovskii L.V.} {\it An estimate for the probabilities of large deviations.  } Mat. Zametki, (1987), {\bf 12},
145 \ - \ 156. Translated from Russian.

\bibitem{Rozovskii1}
{\sc Rozovskii L.V.} {\it  Large deviations probabilities for sums of independent  random variables
with common distribution  from the domain of attraction  of a stable law. }  Journal of Math. Science, Vol. 93, $  N^o 3,$  (1999), 421 \ - \ 433.
 Translation from Theor. Veroyatn. i Primenen., {\bf 42}, 3, (1998),  454 - 481. (in  Russian),

\bibitem{Rozovskii2}
{\sc Rozovskii L.V.} {\it  Large deviations of sums of independent  random variables  from the domain
of attraction  of  non-symmetric  stable law. } (in Russian).  Teor. Veroyatn.  i Primenen., {\bf 42 (3),} (1997), 496 \ - \ 536.

\bibitem{Tsallis1}
{\sc  Tsallis Constantino.} {\it Possible generalization of Boltzmann-Gibbs statistics.} Journal of Statistical
Physics, 1988, {\bf 52,} 479 \ – \ 487.

\bibitem{Varadhan1}
{\sc Varadhan, S.R.S.} (1984) {\it Large Deviations and Applications.} Philadelphia: SIAM.
Mathematical Reviews (MathSciNet): MR758258

\end{thebibliography}
\end{document}